\documentclass[11pt]{amsart}
\usepackage{amscd,amssymb}
\usepackage{amsthm,amsmath,amssymb}
\usepackage[matrix,arrow]{xy}

\newtheorem{theorem}[equation]{Theorem}

\newtheorem{lemma}[equation]{Lemma}

\theoremstyle{definition}
\newtheorem{example}[equation]{Example}

\theoremstyle{remark}
\newtheorem{remark}[equation]{Remark}

\author{Ivan Cheltsov}

\title{On nodal sextic fivefold}

\address{\begin{tabbing}
\hspace*{28 em}\=\kill
Steklov Institute of Mathematics \>School of Mathematics\\
8 Gubkin street, Moscow 117966   \>The University of Edinburgh\\
Russia                           \>Kings Buildings,  Mayfield Road\\
                                 \> Edinburgh EH9 3JZ, UK\\
\texttt{cheltsov@yahoo.com}      \>\\
                                 \>\texttt{I.Cheltsov@ed.ac.uk}
\end{tabbing}}

\begin{document}

\begin{abstract}
We prove the birational superrigidity and nonrationality of a
hypersurface in $\mathbb{P}^{6}$ of degree $6$ having at most
isolated ordinary double points.
\end{abstract}

\maketitle

\section{Introduction.}
\label{section:introduction}

In many cases the only known way to prove the nonrationality of a
Fano variety\footnote{All varieties are assumed to be projective,
normal and defined over $\mathbb{C}$.} is to prove its birational
rigidity\footnote{Let $V$ be a Fano variety with terminal
$\mathbb{Q}$-factorial singularities and
$\mathrm{rk}\,\mathrm{Pic}(V)=1$. Then $V$ is called birationally
rigid if it is not bi\-ra\-ti\-o\-nal to the following varieties:
a variety $Y$ such that there is a morphism $\tau:Y\to Z$ whose
general fiber has negative Kodaira dimension and
$\mathrm{dim}(Y)\ne\mathrm{dim}(Z)\ne 0$; a Fano variety of Picard
rank $1$ having terminal $\mathbb{Q}$-factorial singularities that
is not biregular to $V$. The variety $V$ is called birationally
superrigid if it is birationally rigid and
$\mathrm{Bir}(V)=\mathrm{Aut}(V)$.}. Many counterexamples to the
L\"uroth problem are obtained by proving the birational rigidity
of Fano 3-folds (see \cite{IsMa71}). Moreover, birational rigidity
is the only\footnote{A priori the method of J.Koll\'ar can be
applied to construct explicit examples of nonrational Fano
varieties, but a posteriori there is only one case of such
explicit application (see \cite{Ko96}, \cite{CoKoSm03}).} known
way to prove the nonrationality of an explicitly given Fano
$n$-fold for $n>3$.

Birational rigidity is proved in the following cases:
\begin{itemize}
\item for some smooth Fano 3-folds (see \cite{IsMa71}, \cite{Is80b}, \cite{IsPu96});%
\item for many singular Fano 3-folds (see \cite{Pu88b}, \cite{Pu97}, \cite{Gr98a}, \cite{CPR}, \cite{CoMe02}, \cite{Me03});%
\item for many smooth Fano $n$-folds (see \cite{Pu87}, \cite{Pu98a}, \cite{Pu00b}, \cite{Ch00b}, \cite{Pu01}, \cite{Pu02a}, \cite{Pu03b}, \cite{dFEM03}, \cite{Ch03b}, \cite{Ch04}), $n>3$;%
\item for some singular Fano $n$-folds (see \cite{Pu88b}, \cite{Pu97}, \cite{Pu02c}, \cite{Pu03a}, \cite{Ch04}), $n>3$.%
\end{itemize}

Let $X$ be a hypersurface in $\mathbb{P}^{6}$ of degree $6$ such
that the only singularities of $X$ are isolated ordinary double
points. Then $-K_{X}\sim\mathcal{O}_{\mathbb{P}^{6}}(1)\vert_{X}$,
the variety $X$ is a Fano 5-fold with $\mathbb{Q}$-factorial
terminal singularities and $\mathrm{rk}\,\mathrm{Pic}(X)=1$ (see
\cite{CalLy94}). In this paper we prove the following result.

\begin{theorem}
\label{theorem:main} The hypersurface $X$ is birationally
superrigid.
\end{theorem}

In the smooth case the claim of Theorem~\ref{theorem:main} is
proved in \cite{Ch00b}. In fact, one can use
Theorem~\ref{theorem:main} to construct explicit examples of
nonrational singular hypersurfaces.

\begin{example}
\label{example:single-point} The singularities of the hypersurface
$$
x_{0}^{4}(x_{1}^{2}+x_{2}^{2}+x_{3}^{2}+x_{4}^{2}+x_{5}^{2}+x_{6}^{2})=
x_{1}^{6}+x_{2}^{6}+x_{3}^{6}+x_{4}^{6}+x_{5}^{6}+x_{6}^{6}\subset\mathbb{P}^{6}\cong\mathrm{Proj}\big(\mathbb{C}[x_{0},\ldots,x_{6}]\big).%
$$
consist of a single ordinary double point, which implies that it
is nonrational by Theorem~\ref{theorem:main}.
\end{example}

\begin{example}
\label{example:plane-odd} Let $X$ be a hypersurface
$$
\sum_{i=0}^{2}a_{i}(x_{0},\ldots,x_{6})b_{i}(x_{0},\ldots,x_{6})=0\subset\mathbb{P}^{6}\cong\mathrm{Proj}\big(\mathbb{C}[x_{0},\ldots,x_{6}]\big),%
$$
where $a_{i}$ and $b_{i}$ are general homogeneous polynomials of
degree $3$. Then $X$ has $729$ isolated ordinary double points. In
particular, the hypersurface $X$ is nonrational by
Theorem~\ref{theorem:main}.
\end{example}

It should be pointed out that the claim of
Theorem~\ref{theorem:main} can be considered as a five-dimensional
generalization of the birational rigidity of a
$\mathbb{Q}$-factorial quartic 3-fold having isolated ordinary
double points (see \cite{IsMa71}, \cite{Pu88b}, \cite{Me03}). The
claim of Theorem~\ref{theorem:main} is relevant to \cite{Pu02c}
and \cite{Pu03a}, but one can not use \cite{Pu02c} and
\cite{Pu03a} to produce explicit examples of nonrational Fano
hypersurfaces.

\smallskip

The author is very grateful to I.Aliev, A.Corti, M.Grinenko,
V.Is\-kov\-skikh, J.Park, Yu.Pro\-kho\-rov and V.Sho\-ku\-rov for
fruitful conversations. The author would like to cordially thank
the referee who pointed out the way how to strengthen the original
claim of Lemma~\ref{lemma:6-n-square-dimension-5}, which allowed
to remove a redundant assumption in the original claim of
Theorem~\ref{theorem:main}.

\section{The Noether--Fano--Iskovskikh inequality.}
\label{section:Noether-Fano-Iskovskikh}

Let $X$ be a Fano variety with terminal $\mathbb{Q}$-factorial
singularities  such that $\mathrm{rk}\,\mathrm{Pic}(X)=1$, but the
variety $X$ is not birationally superrigid. Then the following
result holds (see \cite{Co95}).

\begin{theorem}
\label{theorem:Nother-Fano} There is a linear system $\mathcal{M}$
on the variety $X$ whose base locus has codimension at least $2$,
and the singularities of the log pair $(X, \gamma\mathcal{M})$ are
not canonical, where $\gamma$ is a positive rational number such
that the equivalence $K_{X}+\gamma\mathcal{M}\sim_{\mathbb{Q}} 0$
holds.
\end{theorem}

In the rest of the section we prove
Theorem~\ref{theorem:Nother-Fano}. Let $\rho:X\dasharrow Y$ be a
birational map such that the rational map $\rho$ is not biregular
and one of the following holds:
\begin{itemize}
\item the variety $Y$ is a Fano variety with terminal
$\mathbb{Q}$-factorial singularities such that the equality $\mathrm{rk}\,\mathrm{Pic}(Y)=1$ holds (the Fano case);%
\item the variety $Y$ is smooth, and there is a morphism
$\tau:Y\to Z$ whose general fiber has negative Kodaira dimension
and $\mathrm{dim}(Y)\ne\mathrm{dim}(Z)\ne 0$ (the fibration case).
\end{itemize}

Let us consider a commutative diagram
$$
\xymatrix{
&&W\ar@{->}[ld]_{\alpha}\ar@{->}[rd]^{\beta}&&\\%
&X\ar@{-->}[rr]_{\rho}&&Y,&}
$$
such that $W$ is smooth, $\alpha$ and $\beta$ are birational
morphisms.  In the Fano case let $\mathcal{D}$ be the complete
linear system $|-rK_{Y}|$ for $r\gg 0$, in the fibration case let
$\mathcal{D}$ be the complete linear system $|\tau^{*}(H)|$, where
$H$ is a very ample divisor on $Z$. Let $\mathcal{M}$ be a proper
transform on the variety $X$ of the linear system $\mathcal{D}$.
Now choose a positive rational number $\gamma$ such that the
equivalence $K_{X}+\gamma\mathcal{M}\sim_{\mathbb{Q}} 0$ holds.
Suppose that the singularities of the log pair $(X,
\gamma\mathcal{M})$ are not canonical. Let us show that this
assumption leads to a contradiction.

Let $\mathcal{B}$ be a proper transform on $W$ of the linear
system $\mathcal{M}$. Then
$$
\sum_{i=1}^{k}a_{i}F_{i}\sim_{\mathbb{Q}}\alpha^{*}(K_{X}+\gamma\mathcal{M})+\sum_{i=1}^{k}a_{i}F_{i}\sim_{\mathbb{Q}}
K_{W}+\gamma\mathcal{B}\sim_{\mathbb{Q}}\beta^{*}(K_{Y}+\gamma\mathcal{D})+\sum_{i=1}^{l}b_{i}G_{i},%
$$
where $F_{j}$ is a $\beta$-exceptional divisor, $G_{i}$ is an
$\alpha$-exceptional divisor, $a_{i}$ is a nonnegative rational
number, and $b_{i}$ is a positive rational number. Let $n$ be a
sufficiently big and sufficiently divisible natural number. Then
$$
1=h^{0}\Big(\mathcal{O}_{W}\big(\sum_{j=1}^{k}na_{j}F_{j}\big)\Big)=h^{0}\Big(\mathcal{O}_{W}\big(\beta^{*}(nK_{Y}+n\gamma\mathcal{D})+\sum_{i=1}^{l}nb_{i}G_{i}\big)\Big),%
$$
but
$h^{0}(\mathcal{O}_{W}(\beta^{*}(nK_{Y}+\gamma\mathcal{D})+\sum_{i=1}^{l}nb_{i}G_{i}))=0$
in the fibration case. Hence, the fibration case is impossible. In
the Fano case the equality
$h^{0}(\mathcal{O}_{W}(\beta^{*}(nK_{Y}+\gamma\mathcal{D})+\sum_{i=1}^{l}nb_{i}G_{i}))=1$
implies that $\gamma=1/r$. Thus, we have
$$
\sum_{i=1}^{k}a_{i}F_{i}\sim_{\mathbb{Q}}\sum_{i=1}^{l}b_{i}G_{i},%
$$
and it follows from Lemma~2.19 in \cite{Ko91} that
$\sum_{i=1}^{k}a_{i}F_{i}=\sum_{i=1}^{l}b_{i}G_{i}$, which implies
that the singularities of the log pair $(X, \gamma\mathcal{M})$
are terminal.

There is a rational number $\mu>\gamma$ such that both log pairs
$(X, \mu\mathcal{M})$ and $(X, \mu\mathcal{B})$ have terminal
singularities. Hence, we have
$$
\alpha^{*}(K_{X}+\mu\mathcal{M})+\sum_{i=1}^{k}a^{\prime}_{i}F_{i}\sim_{\mathbb{Q}}
K_{W}+\mu\mathcal{B}\sim_{\mathbb{Q}}\beta^{*}(K_{Y}+\mu\mathcal{D})+\sum_{i=1}^{l}b^{\prime}_{i}G_{i},%
$$
where $a^{\prime}_{i}$ and $b^{\prime}_{i}$ are positive rational
numbers. Let $n$ be a sufficiently big and divisible natural
number, and $\psi:W\dasharrow U$ be a map given by the linear
system $|nK_{W}+n\mu\mathcal{B}|$. Then $\psi\circ\beta^{-1}$ is
an isomorphism, because the divisor $n(K_{Y}+\mu\mathcal{D})$ is
very ample, but the divisor $\sum_{i=1}^{l}nb^{\prime}_{i}G_{i}$
is effective and $\beta$-exceptional. Similarly, we get
$\psi\circ\alpha^{-1}$ is an isomorphism. Hence, the birational
map $\rho$ is an isomorphism, which is a contradiction. Thus, we
proved Theorem~\ref{theorem:Nother-Fano}.

\section{The lemma of Corti.}
\label{section:Corti}

Let $X$ be a variety, $O$ be an isolated ordinary double point on
$X$, $B_{X}$ be an effective $\mathbb{Q}$-Cartier divisor on the
variety $X$, $\pi:W\to X$ be a blow up of $O$, $E$ be a
$\pi$-exceptional divisor, $B_{W}$ be a proper transform of the
divisor $B_{X}$ on the variety $W$. Then the equivalence
$$
\pi^{*}(B_{X})\sim_{\mathbb{Q}}B_{W}+\mathrm{mult}_{O}(B_{X})E
$$
holds, where $\mathrm{mult}_{O}(B_{X})$ is a non-negative rational
number. Suppose that $\mathrm{dim}(X)\geqslant 3$ and the
singularities of the log pair $(X, B_{X})$ are not canonical in
the point $O$. Then elementary calculations imply
$\mathrm{mult}_{O}(B_{X})>1/2$. The following result is implied by
Theorem~3.10 in \cite{Co00}.

\begin{lemma}
\label{lemma:Corti} The inequality $\mathrm{mult}_{O}(B_{X})>1$
holds.
\end{lemma}

In the rest of the section we prove Lemma~\ref{lemma:Corti}.
Suppose that $\mathrm{mult}_{O}(B_{X})\leqslant 1$. Let us show
that this assumption leads to a contradiction. Replacing the
divisor $B_{X}$ by $(1-\epsilon)B_{X}$ for some positive
sufficiently small rational $\epsilon$, we may assume that
$\mathrm{mult}_{O}(B_{X})<1$. Moreover, taking sufficiently
general hyperplane sections of $X$, we may assume that
$\mathrm{dim}(X)=3$ due to Theorem~17.6 in \cite{Ko91}.

\begin{lemma}
\label{lemma:quadric-surface} Let $S$ be a surface
$\mathbb{P}^{1}\times \mathbb{P}^{1}$, and $B_{S}$ be an effective
divisor on $S$ of bi-degree$(a, b)$, where $a$ and $b$ are
rational numbers in $[0, 1)$. Then the log pair $(S, B_{S})$ has
log-terminal singularities.
\end{lemma}

\begin{proof}
Suppose that the singularities of $(S, B_{S})$ are not
log-terminal. Then the locus of log ca\-nonical singularities
$\mathrm{LCS}(S, B_{S})$ is not empty and consists of points of
the surface $S$. Hence, the locus $\mathrm{LCS}(S, F+B_{S})$ is
not connected, where $F$ is a sufficiently general fiber of the
projection of the surface $S$ to $\mathbb{P}^{1}$. The later
contradicts Theorem~17.4 in \cite{Ko91}.
\end{proof}

The inequality $\mathrm{mult}_{O}(B_{X})<1$ and the equivalence
$$
K_{W}+B_{W}\sim_\mathbb{Q}
\pi^{*}\big(K_{X}+B_{X}\big)+\big(1-\mathrm{mult}_{O}(B_{X})\big)E,
$$
imply that there is a proper irreducible subvariety $Z\subset E$
such that the log pair $(W, B_{W})$ is not canonical in the
generic point of $Z$. Hence the singularities of the log pair $(E,
B_{W}\vert_{E})$ are not log terminal by Theorem~17.6 in
\cite{Ko91}, which is impossible by
Lemma~\ref{lemma:quadric-surface}.

\section{Main inequalities.}
\label{section:main-inequality}

Let $X$ be a variety, $O$ be an isolated ordinary double point on
$X$, $\mathcal{M}$ be a liner system on the variety $X$ having no
base components, and $r=\mathrm{dim}(X)\geqslant 4$. Let $\pi:V\to
X$ be a blow up of $X$ at the point $O$, $E$ be a
$\pi$-exceptional divisor, and let $\mathcal{B}$ be a proper
transform of the linear system $\mathcal{M}$ on the variety $V$.
Then the divisor $E$ can be identified with a smooth quadric
hypersurface in $\mathbb{P}^{r}$, and the equivalence
$$
\mathcal{B}\sim \pi^{*}(\mathcal{M})-\mathrm{mult}_{O}(\mathcal{M})E%
$$
holds for some natural number $\mathrm{mult}_{O}(\mathcal{M})$. It
should be pointed out that $\mathrm{mult}_{O}(\mathcal{M})$ is
different from the scheme-theoretic multiplicity of a general
surface of $\mathcal{M}$ in the point $O$.

Let $S_{1}$ and $S_{2}$ be general divisors in the linear system
$\mathcal{M}$, and $H_{i}$ be a general hyperplane section of $X$
passing through $O$, where $i=1,\ldots,r-2$. We can define
$\mathrm{mult}_{O}(S_{i})$ and $\mathrm{mult}_{O}(H_{i})$ in the
same way as we defined the number
$\mathrm{mult}_{O}(\mathcal{M})$. Let $\hat{S}_{i}$ and
$\hat{H}_{i}$ be proper transforms on the variety $V$ of the
divisors $S_{i}$ and $H_{i}$ respectively. Then we can put
$$
\mathrm{mult}_{O}\big(S_{1}\cdot S_{2}\big)=2\mathrm{mult}^{2}_{O}(S_{i})+\sum_{P\in E}\mathrm{mult}_{P}(\hat{S}_{1}\cdot\hat{S}_{2})\mathrm{mult}_{P}(\hat{H}_{1})\cdots\mathrm{mult}_{P}(\hat{H}_{r-2}).%
$$

\begin{remark}
\label{remark:remark} The inequality $\mathrm{mult}_{O}(S_{1}\cdot
S_{2})\geqslant
2\mathrm{mult}^{2}_{O}(S_{i})+\mathrm{mult}_{Z}(\hat{S}_{1}\cdot\hat{S}_{2})$
holds for any irreducible subvariety $Z\subset E$ of codimension
one.
\end{remark}

\begin{example}
\label{example:sextic-double-points-I} Let $X$ be a hypersurface
in $\mathbb{P}^{6}$ of degree $6$ such that the singularities of
the hypersurface $X$ consist of a finite number of isolated
ordinary double points, and let $O$ be a singular point of the
variety $X$. Then the groups $\mathrm{Cl}(X)$ and
$\mathrm{Pic}(X)$ are generated by a hyperplane section $H$ of the
hypersurface $X$ (see \cite{CalLy94}), which implies that
$S_{i}\sim nH$ for some natural number $n$. Moreover, the
inequality $\mathrm{mult}_{O}(S_{1}\cdot S_{2})\leqslant 6n^{2}$
holds.
\end{example}

Suppose that the singularities of the log pair $(X,
{\frac{1}{n}}\mathcal{M})$ are not canonical in the point $O$, but
they are canonical in a punctured neighborhood of the point $O$.

\begin{lemma}
\label{lemma:6-n-square-big-dimension} Suppose that
$\mathrm{dim}(X)\geqslant 6$. Then $\mathrm{mult}_{O}(S_{1}\cdot
S_{2})>6n^{2}$.
\end{lemma}

\begin{proof}
We prove the inequality $\mathrm{mult}_{O}(S_{1}\cdot
S_{2})>6n^{2}$ only when $\mathrm{dim}(X)=6$, because the proof in
the case $\mathrm{dim}(X)\geqslant 7$ is similar. So suppose that
$\mathrm{dim}(X)=6$. Then
$$
K_{V}+{\frac{1}{n}}\mathcal{B}\sim_{\mathbb{Q}} \pi^{*}\big(K_{X}+{\frac{1}{n}}\mathcal{M}\big)+\Big(4-{\frac{\mathrm{mult}_{O}(\mathcal{M})}{n}}\Big)E.%
$$

Put $\check{X}=\cap_{i=1}^{3}H_{i}$ and
$\check{\mathcal{M}}=\mathcal{M}\vert_{\check{X}}$. Then $O$ is an
isolated ordinary double point on $\check{X}$, and the
singularities of the log pair $(\check{X},
{\frac{1}{n}}\check{\mathcal{M}})$ are not log canonical in the
point $O$ by Theorem~17.6 of the paper \cite{Ko91}. Let
$\check{\pi}:\check{V}\to\check{X}$ be a blow up of the point $O$,
and $\check{E}$ be an exceptional divisor of the birational
morphism $\check{\pi}$. Then the diagram
$$
\xymatrix{
&\check{V}\ar@{->}[d]_{\check{\pi}}\ar@{^{(}->}[rr]&& V\ar@{->}[d]^{\pi}&\\%
&\check{X}\ar@{^{(}->}[rr]&&X&}
$$
is commutative, where the 3-fold $\check{V}$ is identified with a
proper transform of the subvariety  $\check{X}$ on the variety
$V$. In particular, we have $\check{E}=E\cap\check{V}$. The
generality of $H_{i}$ implies
$$
\mathrm{mult}_{O}(\check{\mathcal{M}})=\mathrm{mult}_{O}(\mathcal{M}),
$$
and we may assume that $\mathrm{mult}_{O}(\mathcal{M})<2n$,
because otherwise $\mathrm{mult}_{O}(S_{1}\cdot S_{2})>6n^{2}$.

Let $\mathcal{B}$ be a proper transform of $\mathcal{M}$ on the
variety $V$, and $\check{\mathcal{B}}$ be a proper transform of
the linear system $\check{\mathcal{M}}$ on the 3-fold $\check{V}$.
Then $\check{\mathcal{B}}=\mathcal{B}\vert_{\check{V}}$ and we
have
$$
K_{V}+{\frac{1}{n}}\mathcal{B}+\Big({\frac{\mathrm{mult}_{O}(\mathcal{M})}{n}}-1\Big)E+\hat{H}_{1}+\hat{H}_{2}+\hat{H}_{3}\sim_{\mathbb{Q}} \pi^{*}\big(K_{X}+{\frac{1}{n}\mathcal{M}}+H_{1}+H_{2}+H_{3}\big)%
$$
and
$$
K_{\check{V}}+{\frac{1}{n}}\check{\mathcal{B}}+\Big({\frac{\mathrm{mult}_{O}(\mathcal{M})}{n}}-1\Big)\check{E}\sim_{\mathbb{Q}} \check{\pi}^{*}\big(K_{\check{X}}+{\frac{1}{n}}\check{\mathcal{M}}\big),%
$$
but $\mathrm{mult}_{O}(\mathcal{M})<2n$ implies the existence of
irreducible subvarieties $\Omega\subsetneq E$ and
$\check{\Omega}\subsetneq\check{E}$ such that the singularities of
the log pair $(V,
{\frac{1}{n}}\mathcal{B}+(\mathrm{mult}_{O}(\mathcal{M})/n-1)E)$
are not log canonical in the generic point of $\Omega$, the
singularities of the log pair $(\check{V},
{\frac{1}{n}}\check{\mathcal{B}}+(\mathrm{mult}_{O}(\mathcal{M})/n-1)\check{E})$
are not log canonical in the generic point of $\check{\Omega}$,
and $\check{\Omega}\subseteq\Omega\cap\check{V}$. We have
$\check{\Omega}=\Omega\cap\check{V}$ when
$\mathrm{dim}(\check{\Omega})>0$, and we may assume that $\Omega$
and $\check{\Omega}$ have the greatest possible dimensions among
all subvarieties having such properties. Applying Theorem 17.4 of
\cite{Ko91} to $(\check{V},
{\frac{1}{n}}\check{\mathcal{B}}+(\mathrm{mult}_{O}(\mathcal{M})/n-1)\check{E})$
and $\check{\pi}$ we see that in the case
$\mathrm{dim}(\check{\Omega})=0$ the locus of log canonical
singularities
$$
\mathrm{LCS}\Big(\check{V},
{\frac{1}{n}}\check{\mathcal{B}}+\big(\mathrm{mult}_{O}(\mathcal{M})/n-1\big)\check{E}\Big)
$$
consists of a single point $\check{\Omega}$ in the neighborhood of
the divisor $\check{E}$. In particular, we have
$\check{\Omega}=\Omega\cap\check{V}$.

Suppose that $\mathrm{dim}(\check{\Omega})=0$. Then
$\check{\Omega}=\Omega\cap\check{V}$ implies that $\Omega$ is a
linear subspace in $\mathbb{P}^{6}$ of codimension $3$ that is
contained in the smooth quadric hypersurface
$E\subset\mathbb{P}^{6}$, which is impossible by the Lefschetz
theorem. Hence, the inequality
$\mathrm{dim}(\check{\Omega})\geqslant 1$ holds, which implies
$\mathrm{dim}(\Omega)=4$.

The singularities of the log pair $(V,
{\frac{1}{n}}\mathcal{B}+(\mathrm{mult}_{O}(\mathcal{M})/n-1)E)$
are not log canonical in the generic point of the subvariety
$\Omega\subset E$ of dimension $4$. Hence, we can apply
Theorem~3.1 of \cite{Co00} to the log pair $(V,
{\frac{1}{n}}\mathcal{B}+(\mathrm{mult}_{O}(\mathcal{M})/n-1)E)$
in the generic point of $\Omega$. The latter gives
$$
\mathrm{mult}_{\Omega}(\hat{S}_{1}\cdot\hat{S}_{2})>4\big(2n^{2}-n\mathrm{mult}_{O}(\mathcal{M})\big),%
$$
where $\hat{S}_{i}$ is a proper transform of $S_{i}$ on the
variety $V$. Hence, the inequalities
$$
\mathrm{mult}_{O}(S_{1}\cdot S_{2})\geqslant
2\mathrm{mult}_{O}(\mathcal{M})^{2}+\mathrm{mult}_{\Omega}(\hat{S}_{1}\cdot\hat{S}_{2})>
6n^{2}+2\Big(n-\mathrm{mult}_{O}(\mathcal{M})\Big)^{2}\geqslant
6n^{2}
$$
hold.
\end{proof}

Let $\Delta$ be an effective divisor on the variety $X$ passing
through the point $O$ and $\hat{\Delta}$ be its proper transform
on the variety $V$. Suppose that the divisor $\Delta$ does not
contain irreducible components of the cycle $S_{1}\cdot S_{2}$,
and the divisor $\hat{\Delta}$ does not contain irreducible
components of the cycle $\hat{S}_{1}\cdot \hat{S}_{2}$. Then we
can put
$$
\mathrm{mult}_{O}(S_{1}\cdot S_{2}\cdot\Delta)=2\mathrm{mult}^{2}_{O}(S_{i})\mathrm{mult}_{O}(\Delta)+\sum_{P\in E}\mathrm{mult}_{P}(\hat{S}_{1}\cdot\hat{S}_{2}\cdot\hat{\Delta})\mathrm{mult}_{P}(\hat{H}_{1})\cdots\mathrm{mult}_{P}(\hat{H}_{r-3}),%
$$
which implies $\mathrm{mult}_{O}(S_{1}\cdot S_{2}\cdot\Delta)=
\mathrm{mult}_{O}(S_{1}\vert_{\Delta}\cdot S_{2}\vert_{\Delta})$
in the case when the point $O$ is an isolated ordinary double
point on the divisor $\Delta$.

\begin{lemma}
\label{lemma:6-n-square-dimension-4} Suppose that
$\mathrm{dim}(X)=4$. Then there is a line $\Lambda\subset
E\subset\mathbb{P}^{4}$ such that the strict inequality
$\mathrm{mult}_{O}(S_{1}\cdot S_{2}\cdot \Delta)>6n^{2}$ holds if
 $\Lambda\subset\hat{\Delta}$, and $O$ is an
ordinary double point on $\Delta$.%
\end{lemma}

\begin{proof}
We have $\mathrm{mult}_{O}(\mathcal{M})>n$ by
Lemma~\ref{lemma:Corti}, but
$$
K_{V}+{\frac{1}{n}}\mathcal{B}\sim_{\mathbb{Q}} \pi^{*}\big(K_{X}+{\frac{1}{n}}\mathcal{M}\big)+\Big(2-{\frac{\mathrm{mult}_{O}(\mathcal{M})}{n}}\Big)E.%
$$

Suppose that $O$ is an ordinary double point on $\Delta$. Put
$\bar{S}_{i}=S_{i}\vert_{\Delta}$ and
$\bar{\mathcal{M}}=\mathcal{M}\vert_{\Delta}$. Then the log pair
$(\Delta, {\frac{1}{n}}\bar{\mathcal{M}})$ is not log canonical in
the point $O$ by Theorem~17.6 in \cite{Ko91}.

Let $\tilde{\pi}:\tilde{\Delta}\to\Delta$ be a blow up of $O$, and
$\tilde{E}$ is a $\bar{\pi}$-exceptional divisor. Then the diagram
$$
\xymatrix{
&\tilde{\Delta}\ar@{->}[d]_{\bar{\pi}}\ar@{^{(}->}[rr]&& V\ar@{->}[d]^{\pi}&\\%
&\Delta\ar@{^{(}->}[rr]&&X&}
$$
is commutative, where we can identify $\tilde{\Delta}$ with
$\hat{\Delta}$, and $\tilde{E}=E\cap\tilde{\Delta}$ can be
considered as a nonsingular quadric hypersurface in
$\mathbb{P}^{3}$. The inequality
$\mathrm{mult}_{O}(\bar{\mathcal{M}})\geqslant 2n$ gives
$$
\mathrm{mult}_{O}\big(S_{1}\cdot S_{2}\cdot\Delta\big)=\mathrm{mult}_{O}\big(\bar{S}_{1}\cdot\bar{S}_{2}\big)\geqslant 8n^{2},%
$$
hence, we may assume that
$\mathrm{mult}_{O}(\bar{\mathcal{M}})<2n$.

Let $\tilde{\mathcal{M}}$ be a proper transform of the linear
system $\bar{\mathcal{M}}$ on $\tilde{\Delta}$. Then
$\mathrm{mult}_{O}(\bar{\mathcal{M}})<2n$ implies the existence of
an irreducible subvariety $\Xi\subsetneq\tilde{E}$ such that the
singularities of the log pair
$$
\Big(\tilde{\Delta}, {\frac{1}{n}}\tilde{\mathcal{M}}+\big(\mathrm{mult}_{O}(\bar{\mathcal{M}})/n-1\big)\tilde{E}\Big).%
$$
are not log canonical in the generic point of $\Xi$.

Suppose that $\Xi$ is a curve. Let $\tilde{S}_{i}$ be a proper
transform of $\bar{S}_{i}$ on $\tilde{\Delta}$. Then
$$
\mathrm{mult}_{O}\big(\bar{S}_{1}\cdot \bar{S}_{2}\big)\geqslant 2\mathrm{mult}_{O}\big({\mathcal{M}}\big)^{2}+\mathrm{mult}_{\Xi}\big(\tilde{S}_{1}\cdot\tilde{S}_{2}\big),%
$$
but Theorem~3.1 of \cite{Co00} applied to the log pair
$(\tilde{\Delta},
{\frac{1}{n}}\tilde{\mathcal{M}}+(\mathrm{mult}_{O}(\bar{\mathcal{M}})/n-1)\tilde{E})$
in the generic point of $\Xi$ implies that the inequality
$$
\mathrm{mult}_{\Xi}\big(\tilde{S}_{1}\cdot\tilde{S}_{2}\big)>4\big(2n^{2}-n\mathrm{mult}_{O}(\bar{\mathcal{M}})\big)
$$
holds. Hence, the inequalities
$$
\mathrm{mult}_{O}\big(\bar{S}_{1}\cdot\bar{S}_{2}\big)>2\mathrm{mult}_{O}^{2}(\bar{\mathcal{M}})+4\big(2n^{2}-n\mathrm{mult}_{O}(\bar{\mathcal{M}})\big)\geqslant 6n^{2}%
$$
hold. Thus, we may assume that $\Xi$ is a point.

Suppose that the divisor $\Delta$ is a sufficiently general
hyperplane section of $X$ passing through the point $O$. Then
applying Theorem 17.4 of \cite{Ko91} to the log pair
$(\tilde{\Delta},
{\frac{1}{n}}\tilde{\mathcal{M}}+(\mathrm{mult}_{O}(\bar{\mathcal{M}})/n-1)\tilde{E})$
and the morphism $\tilde{\pi}$ we see that one of the following
holds:
\begin{itemize}
\item the singularities of the log pair $(V,
{\frac{1}{n}}\mathcal{B}+(\mathrm{mult}_{O}(\mathcal{M})/n-1)E)$
are not log canonical in the generic point of some surface that is contained in the divisor $E$;%
\item there is a line $\Lambda\subset E\subset\mathbb{P}^{4}$ such
that the singularities of $(V,
{\frac{1}{n}}\mathcal{B}+(\mathrm{mult}_{O}(\mathcal{M})/n-1)E)$
are not log canonical in the generic point of line $\Lambda$ and
$\Xi=\Lambda\cap \hat{\Delta}$.%
\end{itemize}

In the case when the singularities of the log pair $(V,
{\frac{1}{n}}\mathcal{B}+(\mathrm{mult}_{O}(\mathcal{M})/n-1)E)$
are not log canonical in the generic point of some surface
contained in $E$, the previous arguments implies the inequality
$\mathrm{mult}_{O}(\bar{S}_{1}\cdot\bar{S}_{2})>6n^{2}$. Thus, we
may assume that there is a line $\Lambda\subset
E\subset\mathbb{P}^{4}$ such that $\Xi=\Lambda\cap\tilde{\Delta}$
and the singularities of the log pair $(V,
{\frac{1}{n}}\mathcal{B}+(\mathrm{mult}_{O}(\mathcal{M})/n-1)E)$
are not log canonical in the generic point of the curve $\Lambda$.
It should be pointed out that the line $\Lambda$ does not depend
on the choice of the divisor $\Delta$. Therefore, we may assume
that the divisor $\Delta$ is chosen under the additional
assumption $\Lambda\subset\hat{\Delta}$, where we identified
$\hat{\Delta}$ with $\tilde{\Delta}$.

The singularities of the log pair $(\tilde{\Delta},
{\frac{1}{n}}\tilde{\mathcal{M}}+(\mathrm{mult}_{O}(\bar{\mathcal{M}})/n-1)\tilde{E})$
are not log canonical in the generic point of $\Lambda$ by
Theorem~17.6 in \cite{Ko91}, because the boundary
${\frac{1}{n}}\mathcal{B}+(\mathrm{mult}_{O}(\mathcal{M})/n-1)E$
is effective due to the inequality
$\mathrm{mult}_{O}(\mathcal{M})>n$. Hence, we can apply
Theorem~3.1 of \cite{Co00} to the log pair $(\tilde{\Delta},
{\frac{1}{n}}\tilde{\mathcal{M}}+(\mathrm{mult}_{O}(\bar{\mathcal{M}})/n-1)\tilde{E})$
in the generic point of $\Lambda$ to obtain the inequalities
$$
\mathrm{mult}_{O}\big(\bar{S}_{1}\cdot\bar{S}_{2}\big)>2\mathrm{mult}_{O}^{2}(\bar{\mathcal{M}})+4\big(2n^{2}-n\mathrm{mult}_{O}(\bar{\mathcal{M}})\big)\geqslant 6n^{2},%
$$
which conclude the proof.
\end{proof}

Finally, let us prove the following result.

\begin{lemma}
\label{lemma:6-n-square-dimension-5} Suppose that
$\mathrm{dim}(X)=5$ Then
$\mathrm{mult}_{O}(S_{1}\cdot S_{2})>6n^{2}$.%
\end{lemma}

\begin{proof}
Put $\check{X}=H_{1}\cap H_{2}$ and
$\check{\mathcal{M}}=\mathcal{M}\vert_{\check{X}}$. Then $O$ is an
isolated ordinary double point on $\check{X}$, and the
singularities of the log pair $(\check{X},
{\frac{1}{n}}\check{\mathcal{M}})$ are not log canonical in the
point $O$ by Theorem~17.6 of the paper \cite{Ko91}. Let
$\check{\pi}:\check{V}\to\check{X}$ be a blow up of the point $O$,
and $\check{E}$ be an exceptional divisor of the morphism
$\check{\pi}$. Then we can identified $\check{V}$ with a proper
transform of $\check{X}$ on $V$. We have
$$
\mathrm{mult}_{O}\big(S_{1}\cdot S_{2}\big)\geqslant
2\mathrm{mult}^{2}_{O}(\mathcal{M})>6n^{2}
$$
in the case when $\mathrm{mult}_{O}(\mathcal{M})\geqslant 2n$.
Hence, we may assume that $\mathrm{mult}_{O}(\mathcal{M})<2n$.

Let $\check{\mathcal{B}}$ be a proper transform of the linear
system $\check{\mathcal{M}}$ on $\check{V}$. Then
$\check{\mathcal{B}}=\mathcal{B}\vert_{\check{V}}$ and we have
$$
K_{V}+{\frac{1}{n}}\mathcal{B}+\Big({\frac{\mathrm{mult}_{O}(\mathcal{M})}{n}}-1\Big)E+\hat{H}_{1}+\hat{H}_{2}\sim_{\mathbb{Q}} \pi^{*}\big(K_{X}+{\frac{1}{n}\mathcal{M}}+H_{1}+H_{2}\big),%
$$
but
$K_{\check{V}}+{\frac{1}{n}}\check{\mathcal{B}}+(\mathrm{mult}_{O}(\mathcal{M})/n-1)\check{E}\sim_{\mathbb{Q}}
\check{\pi}^{*}(K_{\check{X}}+{\frac{1}{n}}\check{\mathcal{M}})$.
Therefore, there are proper irreducible subvarieties
$\Omega\subsetneq E$ and $\check{\Omega}\subsetneq\check{E}$ such
that $\check{\Omega}\subseteq\Omega\cap\check{V}$ and the
following holds:
\begin{itemize}
\item the log pair $(V,
{\frac{1}{n}}\mathcal{B}+(\mathrm{mult}_{O}(\mathcal{M})/n-1)E)$
is not log canonical in $\Omega$;%
\item the log pair $(\check{V},
{\frac{1}{n}}\check{\mathcal{B}}+(\mathrm{mult}_{O}(\mathcal{M})/n-1)\check{E})$
is not log canonical in $\check{\Omega}$.%
\end{itemize}

We may assume that $\Omega$ and $\check{\Omega}$ have the greatest
possible dimensions among all subvarieties having such properties.
Therefore, we have $\check{\Omega}=\Omega\cap\check{V}$ in the
case when $\mathrm{dim}(\check{\Omega})\geqslant 1$.

Suppose that $\mathrm{dim}(\check{\Omega})\geqslant 1$ holds. Then
$\mathrm{dim}(\Omega)=3$ and we can apply Theorem~3.1 of
\cite{Co00} to the log pair $(V,
{\frac{1}{n}}\mathcal{B}+(\mathrm{mult}_{O}(\mathcal{M})/n-1)E)$
in the generic point of $\Omega$. Therefore, we have
$$
\mathrm{mult}_{\Omega}\big(\hat{S}_{1}\cdot\hat{S}_{2}\big)>4\big(2n^{2}-n\mathrm{mult}_{O}(\mathcal{M})\big),%
$$
which implies that the inequalities
$$
\mathrm{mult}_{O}\big(S_{1}\cdot S_{2}\big)\geqslant 2\mathrm{mult}_{O}^{2}(\mathcal{M})+\mathrm{mult}_{\Omega}\big(\hat{S}_{1}\cdot\hat{S}_{2}\big)>6n^{2}%
$$
hold. Therefore, we may assume that
$\mathrm{dim}(\check{\Omega})=0$.

Applying Theorem 17.4 of \cite{Ko91} to $(\check{V},
{\frac{1}{n}}\check{\mathcal{B}}+(\mathrm{mult}_{O}(\mathcal{M})/n-1)\check{E})$
and $\check{\pi}$ we see that the locus
$$
\mathrm{LCS}\Big(\check{V},
{\frac{1}{n}}\check{\mathcal{B}}+\big(\mathrm{mult}_{O}(\mathcal{M})/n-1\big)\check{E}\Big)
$$
consists of a single point $\check{\Omega}$ in the neighborhood of
the divisor $\check{E}$. Hence, the subvariety $\Omega$ is a plane
in $\mathbb{P}^{5}$. In fact, the subvariety $\Omega$ can not be a
plane\footnote{The referee pointed out to the author that the
subvariety $\Omega$ can not be a plane. We follow the arguments of
the referee to conclude the proof of
Lemma~\ref{lemma:6-n-square-dimension-5}.}. Let us prove the
latter by using the arguments of the original proof of
Lemma~\ref{lemma:Corti} (see Theorem~3.10 in \cite{Co00}).

Let $\breve{X}$ be a general hyperplane section of $X$ passing
through $O$ that is locally given as
$$
xy+zt=0\subset\mathbb{C}^{5}\cong\mathrm{Spec}\big(\mathbb{C}[x,y,z,t,u]\big)
$$
in the neighborhood of $O$, which is given by $x=y=z=t=u=0$. Then
$\breve{X}$ has non-isolated singularities, but we can apply the
previous arguments to $\breve{X}$. Namely, let $\breve{V}$ be a
proper transform of the variety $\breve{X}$ on $V$, and
$\breve{\pi}:\breve{V}\to\breve{X}$ be the induced birational
morphism. Then
$$
K_{\breve{V}}+{\frac{1}{n}}\breve{\mathcal{B}}+\big(\mathrm{mult}_{O}(\mathcal{M})/n-2\big)\breve{E}\sim_{\mathbb{Q}} \breve{\pi}^{*}\big(K_{\breve{X}}+\frac{1}{n}\mathcal{M}\vert_{\breve{X}}\big),%
$$
where $\breve{\mathcal{B}}=\mathcal{B}\vert_{\breve{V}}$, and
$\breve{E}$ is the exceptional divisor of $\breve{\pi}$, which is
a cone over $\mathbb{P}^{1}\times\mathbb{P}^{1}$.

Let $\breve{S}_{x}$ and $\breve{S}_{y}$ be irreducible reduced
Weil divisors on the variety $\breve{X}$ that are given by the
equations $x=t=0$ and $y=t=0$ respectively. Then $\breve{S}_{x}$
and $\breve{S}_{y}$ are not $\mathbb{Q}$-Cartier divisors, but the
divisor $\breve{S}_{x}+\breve{S}_{y}$ is Cartier and given by the
equation $t=0$. Moreover, the equivalence
$$
K_{\breve{V}}+{\frac{1}{n}}\breve{\mathcal{B}}+\big(\mathrm{mult}_{O}(\mathcal{M})/n-1\big)\breve{E}+\breve{H}_{x}+\breve{H}_{y}\sim_{\mathbb{Q}} \breve{\pi}^{*}\big(K_{\breve{X}}+\frac{1}{n}\mathcal{M}\vert_{\breve{X}}+\breve{S}_{x}+\breve{S}_{y}\big),%
$$
holds, where $\breve{H}_{x}$ and $\breve{H}_{y}$ are proper
transforms of $\breve{S}_{x}$ and $\breve{S}_{y}$ on the variety
$\breve{V}$. Then
$$
\mathrm{LCS}\Big(\breve{V}, {\frac{1}{n}}\breve{\mathcal{B}}+\big(\mathrm{mult}_{O}(\mathcal{M})/n-1\big)\breve{E}\Big)=\breve{\Omega},%
$$
where $\breve{\Omega}=\Omega\vert_{\breve{V}}$, because we can
apply the previous arguments to $(\breve{X},
\frac{1}{n}\mathcal{M}\vert_{\breve{X}}+\breve{S}_{x}+\breve{S}_{y})$
due to the generality in the choice of $\breve{X}$. Note, that
$\breve{\Omega}$ is a line on the quadric cone
$\breve{E}\subset\mathbb{P}^{4}$.

There are natural ways to desingularize $\breve{X}$ and
$\breve{V}$. Indeed, consider a commutative diagram
$$
\xymatrix{
&& \breve{U}_{x}\ar@{->}[lld]_{\breve{\gamma}_{x}}\ar@{->}[rrrrrr]^{\breve{\eta}_{x}}&& &&  &&\breve{W}_{x}\ar@{->}[rrd]^{\breve{\alpha}_{x}}&&\\%
\breve{V}&& &&\breve{U}\ar@{->}[rr]^{\breve{\xi}}\ar@{->}[llll]_{\breve{\psi}}\ar@{->}[ull]_{\breve{\delta}_{x}}\ar@{->}[dll]^{\breve{\delta}_{y}}&&\breve{W}\ar@{->}[rrrr]^{\breve{\phi}}\ar@{->}[urr]^{\breve{\beta}_{x}}\ar@{->}[drr]_{\breve{\beta}_{y}}&& &&\breve{X},\\%
&& \breve{U}_{y}\ar@{->}[llu]^{\breve{\gamma}_{y}}\ar@{->}[rrrrrr]_{\breve{\eta}_{y}} && && &&\breve{W}_{y}\ar@{->}[rru]_{\breve{\alpha}_{y}}&&}%
$$
where we have the following notations:
\begin{itemize}
\item $\breve{\phi}$ is a blow up of the ideal sheaf of the curve $x=y=z=t=0$;%
\item $\breve{\alpha}_{x}$ and $\breve{\alpha}_{y}$ are blow ups of the ideal sheaves of $\breve{S}_{x}$ and $\breve{S}_{y}$ respectively;%
\item $\breve{\beta}_{x}$ and $\breve{\beta}_{y}$ are blow ups of the exceptional surfaces of $\breve{\alpha}_{x}$ and $\breve{\alpha}_{y}$ respectively;%
\item $\breve{\xi}$, $\breve{\beta}_{x}$, $\breve{\beta}_{y}$ are blow ups of the fibers of $\phi$, $\breve{\alpha}_{x}$, $\breve{\alpha}_{y}$ over the point $O$ respectively;%
\item $\breve{\psi}$ is a blow up of the ideal sheaf of the proper transform of $x=y=z=t=0$;%
\item $\breve{\gamma}_{x}$ and $\breve{\gamma}_{y}$ are blow ups of the ideal sheaves of $\breve{H}_{x}$ and $\breve{H}_{y}$ respectively;%
\item $\breve{\delta}_{x}$ and $\breve{\delta}_{y}$ are blow ups of the exceptional surfaces of $\breve{\gamma}_{x}$ and $\breve{\gamma}_{y}$ respectively.%
\end{itemize}

The varieties $\breve{W}$, $\breve{W}_{x}$, $\breve{W}_{y}$,
$\breve{U}$, $\breve{U}_{x}$, $\breve{U}_{y}$ are smooth by
construction. Moreover, the birational morphisms
$\breve{\alpha}_{x}$, $\breve{\alpha}_{y}$, $\breve{\gamma}_{x}$,
$\breve{\gamma}_{y}$ are small\footnote{A birational morphism is
called small if it does not contract any divisor.}, and
$\breve{\pi}\circ\breve{\psi}=\breve{\phi}\circ\breve{\xi}$. Let
$\breve{F}$ be the $\breve{\xi}$-exceptional divisor. Then
$$
\breve{F}\cong\mathbb{P}\big(\mathcal{O}_{\mathbb{P}^{1}\times\mathbb{P}^{1}}\oplus\mathcal{O}_{\mathbb{P}^{1}\times\mathbb{P}^{1}}(1)\big),
$$
where $\mathcal{O}_{\mathbb{P}^{1}\times\mathbb{P}^{1}}(1)$ is a
hyperplane section of the quadric
$\mathbb{P}^{1}\times\mathbb{P}^{1}$ with respect to the natural
embedding into $\mathbb{P}^{3}$. The induced morphism
$\breve{\xi}\vert_{\breve{F}}$ is the natural projection to
$\mathbb{P}^{1}\times\mathbb{P}^{1}$, the induced morphisms
$\breve{\eta}_{x}\circ\breve{\delta}_{x}\vert_{\breve{F}}$ and
$\breve{\eta}_{y}\circ\breve{\delta}_{y}\vert_{\breve{F}}$ are
projections to $\mathbb{P}^{1}$, the morphisms
$\breve{\delta}_{x}\vert_{\breve{F}}$ and
$\breve{\delta}_{y}\vert_{\breve{F}}$ are contractions of the
exceptional section of $\breve{F}$ to curves, and
$\breve{\psi}\vert_{\breve{F}}$ is the con\-trac\-tion of the
exceptional section of the surface $\breve{F}$ to the vertex of
the cone $\breve{E}$, where $\breve{E}=\breve{\psi}(\breve{F})$.

The subvariety $\breve{\Omega}$ is a line on the quadric cone
$\breve{E}\subset\mathbb{P}^{4}$ that does not pass through the
vertex of the quadric cone $\breve{E}$, but
$(\breve{H}_{x}+\breve{H}_{y})\cdot\breve{\Omega}=1$. We may
assume that $\breve{H}_{x}\cdot\breve{\Omega}=0$ and
$\breve{H}_{y}\cdot\breve{\Omega}=1$.

Let $\breve{D}_{x}$ and $\breve{D}_{y}$ be the proper transforms
of $\breve{H}_{x}$ and $\breve{H}_{y}$ on $\breve{U}_{y}$
respectively, and $\breve{\Gamma}$ be the proper transform of
$\breve{\Omega}$ on the variety $\breve{U}_{y}$. Then
$\breve{D}_{x}\cdot\breve{\Gamma}=0$ and
$\breve{D}_{y}\cdot\breve{\Gamma}=1$. Moreover, we have
$$
K_{\breve{U}_{y}}+{\frac{1}{n}}\breve{\mathcal{D}}+\big(\mathrm{mult}_{O}(\mathcal{M})/n-1\big)\breve{G}+\breve{D}_{x}+\breve{D}_{y}
\sim_{\mathbb{Q}}\big(\breve{\pi}\circ\breve{\gamma}_{y}\big)^{*}\Big(K_{\breve{X}}+\frac{1}{n}\mathcal{M}\vert_{\breve{X}}+\breve{S}_{x}+\breve{S}_{y}\Big),%
$$
where $\breve{\mathcal{D}}$ and $\breve{G}$ are proper transforms
of the linear system $\breve{\mathcal{B}}$ and exceptional divisor
$\breve{E}$ on the variety $\breve{U}_{y}$. The morphism
$\breve{\eta}_{y}$ contracts the divisor $\breve{G}$, but the
morphism $\breve{\eta}_{y}\vert_{\breve{G}}$ is a
$\mathbb{P}^{2}$-bundle.

Let $\breve{Y}$ be a general fiber of
$\breve{\eta}_{y}\vert_{\breve{G}}$. Then
$\breve{Y}\cap\breve{D}_{x}$ is a line in
$\breve{Y}\cong\mathbb{P}^{2}$, the intersection
$\breve{\Gamma}\cap\breve{Y}$ is a point that is not contained in
$\breve{Y}\cap\breve{D}_{x}$, and
$\breve{Y}\cap\breve{D}_{y}=\varnothing$. Therefore, in the
neighborhood of the fiber $Y$ of the morphism $\breve{\eta}_{y}$
the locus of log canonical singularities
$$
\mathrm{LCS}\Big(\breve{U}_{y},
{\frac{1}{n}}\breve{\mathcal{D}}+\big(\mathrm{mult}_{O}(\mathcal{M})/n-1\big)\breve{G}+\breve{D}_{x}+\breve{D}_{y}\Big)
$$
consists of $\breve{\Gamma}$ and $\breve{D}_{x}$, which
contradicts Theorem~17.4 in \cite{Ko91}, because
$\breve{\Gamma}\cap \breve{D}_{x}=\varnothing$.
\end{proof}

\section{The proof of Theorem~\ref{theorem:main}.}
\label{section:birational-rigidity}

Let $X$ be a hypersurface in $\mathbb{P}^{6}$ of degree $6$ having
at most isolated ordinary double points, which is not birationally
superrigid. Let us show that this assumption leads to a
contradiction.

It follows from Theorem~\ref{theorem:Nother-Fano} that there is a
linear system $\mathcal{M}$ on the hypersurface $X$ that does not
have fixed components such that the singularities of the log pair
$(X, {\frac{1}{m}}\mathcal{M})$ are not canonical, where $m$ is a
natural number such that the rational equivalence $\mathcal{M}\sim
-mK_{X}$ holds.

Let $Z$ be a proper irreducible subvariety of $X$ such that the
log pair $(X, {\frac{1}{m}}\mathcal{M})$ is not canonical in the
generic point of $Z$, and $Z$ has maximal dimension among the
subvarieties of $X$ with such property. Then
$\mathrm{dim}(Z)\leqslant 1$ by Theorem~2 in \cite{Pu95}.

Suppose that either $\mathrm{dim}(Z)\ne 0$ or $Z$ is a smooth
point of the hypersurface $X$. Let $P$ be a sufficiently general
point of $Z$, $V$ be a sufficiently general hyperplane section of
$X$ passing through the point $P$, and
$\mathcal{B}=\mathcal{M}\vert_{V}$. Then $V$ is a smooth
hypersurface in $\mathbb{P}^{5}$ of degree $6$, and the
singularities of $(V, {\frac{1}{m}}\mathcal{B})$ are not canonical
in $P$ by Theorem 17.6 of \cite{Ko91}. Let $S_{1}$ and $S_{2}$ be
sufficiently general divisors in $\mathcal{B}$, and $F=S_{1}\cdot
S_{2}$. Then
$$
\mathrm{dim}\big\{O\in F\vert\
\mathrm{mult}_{O}(F)>m\big\}\leqslant 1
$$
by Proposition~5 in \cite{Pu02a}. Let $Y$ be a sufficiently
general hyperplane section of $V$ passing through the point $P$,
and $\mathcal{P}=\mathcal{B}\vert_{Y}$. Then $Y$ is a smooth
hypersurface in $\mathbb{P}^{4}$ of degree $6$, and
\begin{equation}
\label{equation:multiplicities-of-hyperplane-cut}
\mathrm{dim}\big\{O\in F\cap Y\vert\ \mathrm{mult}_{O}\big(F\vert_{Y}\big)>m\big\}\leqslant 0%
\end{equation}
by Proposition~4.5 in \cite{dFEM03}. On the other hand, the
singularities of the log pair $(Y, {\frac{1}{m}}\mathcal{P})$ are
not log canonical in $P$ by Theorem 17.6 of \cite{Ko91}. Let
$\eta:\mathbb{P}^{4}\dasharrow\mathbb{P}^{2}$ be a general
projection. Then
$$
\eta(P)\in\mathrm{LCS}\Big(\mathbb{P}^{2}, {\frac{1}{4m^{2}}}\eta_{*}\big[F\vert_{Y}\big]\Big)%
$$
by Theorem~1.1 in \cite{dFEM03}. Moreover, it follows from the
inequality~\ref{equation:multiplicities-of-hyperplane-cut} and
Proposition~4.7 in \cite{dFEM03} that the singularities of the log
pair $(\mathbb{P}^{2},
{\frac{1}{4m^{2}}}\eta_{*}\big[F\vert_{Y}\big])$ are log terminal
in a punctured neighborhood of the point $\eta(P)$. Hence, the
locus $\mathrm{LCS}(\mathbb{P}^{2},
L+{\frac{1}{4m^{2}}}\eta_{*}[F\vert_{Y}])$ is not connected for a
sufficiently general line $L\subset\mathbb{P}^{2}$, which is
impossible by Theorem 17.4 of \cite{Ko91}, because
$$
K_{\mathbb{P}^{2}}+L+\frac{1}{4m^{2}}\eta_{*}[F\vert_{Y}]\sim_{\mathbb{Q}}-{\frac{1}{2}}L.%
$$

Therefore, we proved that $Z$ is a singular point of $X$. Let
$\pi:U\to X$ be a blow up of the point $Z$, and $E$ be a
$\pi$-exceptional divisor. Then $\mathrm{mult}_{Z}(\mathcal{M})>m$
by Lemma~\ref{lemma:Corti}, but
$$
K_{U}+{\frac{1}{m}}\mathcal{H}\sim_{\mathbb{Q}} \pi^{*}\big(K_{X}+{\frac{1}{m}}\mathcal{M}\big)+\Big(3-{\frac{1}{m}}\mathrm{mult}_{Z}(\mathcal{M})\Big)E,%
$$
where $\mathcal{H}$ is a proper transform of $\mathcal{M}$ on $U$.
Let $M_{1}$ and $M_{2}$ be sufficiently general divisors in the
linear system $\mathcal{M}$. Then the inequality
$$
\mathrm{mult}_{Z}\big(M_{1}\cdot M_{2}\big)>6m^{2}
$$
holds by Lemma~\ref{lemma:6-n-square-dimension-5}. Hence, we have
$$
6m^{2}=M_{1}\cdot M_{2}\cdot H_{1}\cdot H_{2}\cdot H_{3}\geqslant \mathrm{mult}_{Z}\big(M_{1}\cdot M_{2}\big)>6m^{2},%
$$
where $H_{i}$ is a sufficiently general hyperplane section of the
hypersurface $X$ that passes through the point $Z$. The obtained
contradiction proves Theorem~\ref{theorem:main}.

\end{document}